\documentstyle[amssymb]{amsart}

\title{Lebesgue numbers and Atsuji spaces\\in subsystems of second 
order arithmetic}
\author{Mariagnese Giusto}
\address{via Loreto Vecchia 9/10/A\\17100 Savona\\Italy}
\email{giusto@@dm.unito.it}
\author{Alberto Marcone}
\address{Dip.\ di Matematica\\Universit\`a di Torino\\via Carlo 
Alberto 10\\10123 Torino\\Italy}
\email{marcone@@dm.unito.it}

\newtheorem{theorem}{Theorem}[section]
\newtheorem{lemma}[theorem]{Lemma}

\newtheorem{sublemma}{Sublemma}[theorem]

\theoremstyle{definition}
\newtheorem{definition}[theorem]{Definition}
\newtheorem{remark}[theorem]{Remark}

\renewcommand{\P}[2]{$\Pi_{#2}^{#1}$}
\renewcommand{\S}[2]{$\Sigma_{#2}^{#1}$}
\newcommand{\D}[2]{$\Delta_{#2}^{#1}$}
\newcommand{\system}[1]%
{\mbox{\fontfamily{cmss}\fontshape{n}\fontseries{m}%
\selectfont#1}}
\newcommand{\RCA}{\system{RCA}$_0$}
\newcommand{\WKL}{\system{WKL}$_0$}
\newcommand{\ACA}{\system{ACA}$_0$}

\newcommand{\PCA}{\system{\P11-CA}$_0$}

\newcommand{\set}[2]{\{\,{#1}\mid{#2}\,\}}

\newcommand{\cl}[1]{\overline{#1}}
\newcommand{\imply}{\;\longrightarrow\;}
\newcommand{\sse}{\longleftrightarrow}

\renewcommand{\)}{\right)}
\newcommand{\<}{\left\langle}
\renewcommand{\>}{\right\rangle}
\newcommand{\N}{\Bbb N}
\newcommand{\Q}{\Bbb Q}
\newcommand{\R}{\Bbb R}
\newcommand{\A}{\widehat A}
\newcommand{\B}{\widehat B}
\newcommand{\U}{\frak U}
\newcommand{\V}{\frak V}
\newcommand{\seq}{\N^{<\N}}
\newcommand{\Bai}{\N^{\N}}
\newcommand{\Seq}{2^{<\N}}
\newcommand{\Can}{2^{\N}}
\newcommand{\lh}{\operatorname{lh}}
\newcommand{\conc}{{}^\smallfrown}
\newcommand{\f}[1]{{\footnotesize \ref{#1}}}

\begin{document}
\pagestyle{plain}
\maketitle
\begin{abstract}
We study Lebesgue and Atsuji spaces within subsystems of second 
order arithmetic. The former spaces are those such that every open 
covering has a Lebesgue number, while the latter are those such 
that every continuous function defined on them is uniformly 
continuous. The main results we obtain are the following: the 
statement ``every compact space is Lebesgue'' is equivalent to 
\WKL; the statements ``every perfect Lebesgue space is compact'' 
and ``every perfect Atsuji space is compact'' are equivalent to 
\ACA; the statement ``every Lebesgue space is Atsuji'' is provable 
in \RCA; the statement ``every Atsuji space is Lebesgue'' is 
provable in \ACA. We also prove that the statement ``the distance 
from a closed set is a continuous function'' is equivalent to 
\PCA.
\end{abstract}

\section{Introduction}
This paper is part of the program started by Harvey Friedman and 
Steve Simpson and known as {\em reverse mathematics\/}: the aim of 
this program is to understand the role of set existence axioms in 
the development of ordinary mathematics and its present stage 
consists of establishing the weakest subsystem of second order 
arithmetic in which a theorem of ordinary mathematics can be 
proved. The basic reference for this program is Simpson's 
monograph (\cite{sosoa}) while an overview can be found in 
\cite{siZ2}.

We are interested in the theory of complete separable metric 
spaces (also called Polish spaces) and therefore we need to 
develop this theory within weak subsystems of second order 
arithmetic. Such a development is possible by using an appropriate 
coding of these spaces, their subsets and the continuous functions 
among them. This coding is now standard (see e.g.\ 
\cite{bs1,br-thesis,br,st,bs2,sosoa}) and has been used to study 
--- among other things --- various notions of open, closed and 
compact sets, properties of continuous functions and of sequential 
convergence, and basic theorems such as the category theorem, 
various fixed point theorems and the Hahn-Banach theorem for 
separable Banach spaces.

In this paper we concentrate on properties regarding Lebesgue 
numbers of open coverings: these are introduced in basic topology 
(see e.g.\ \cite{kel}). They are a tool for proving that every 
continuous function with compact domain is uniformly continuous 
(see e.g.\ \cite{hy}): by combining lemma \ref{comp->leb} and 
theorem \ref{leb->ats} in this paper we obtain that this proof can 
be carried out in the subsystem \WKL. Lebesgue numbers are also 
used to prove more advanced results in various areas of geometry: 
see e.g.\ \cite{doc} and \cite{kos}.

We call Lebesgue spaces the spaces such that every open covering 
has a Lebesgue number: they turn out to be the spaces $X$ such 
that every continuous function on $X$ is uniformly continuous. 
Spaces having the latter property have been studied for their own 
sake (see \cite{ats,bee1,bee2} and the references quoted therein) 
and are usually called Atsuji spaces.

Let us notice that since every Lebesgue (or Atsuji) metric space 
is complete the restriction imposed by the expressive power of the 
language of second order arithmetic on the spaces we study 
consists solely of forsaking non separable spaces.

The following picture shows the results we obtain by indicating 
the systems needed to prove the various implications between the 
notions of Atsuji and Lebesgue space and these and two notions of 
compactness. The numbers refer to the lemmas or theorems where the 
results are established and the question marks appear beside the 
two implications where we do not know whether our result is 
optimal.

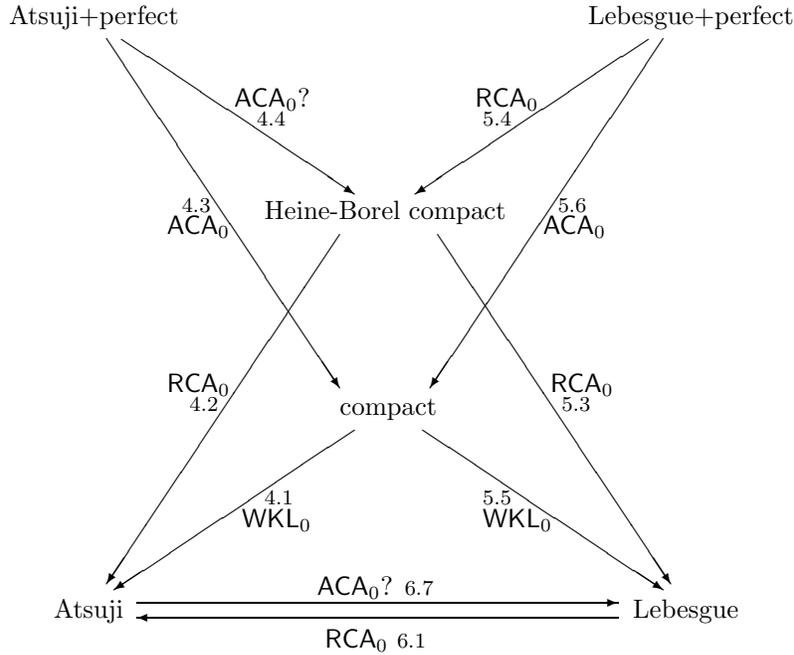
\begin{center}
\setlength{\unitlength}{1cm}
\begin{picture}(10.4,9.1)(0,-8.6)
\put(0,0){Atsuji+perfect} \put(7.7,0){Lebesgue+perfect}
\put(3.4,-2.6){Heine-Borel compact} \put(4.4,-5.2){compact}
\put(0.6,-7.9){Atsuji} \put(8.3,-7.9){Lebesgue}
\put(1.5,-0.2){\vector(3,-2){3.1}} 
\put(8.5,-0.2){\vector(-3,-2){3.1}}
\put(1.3,-0.2){\vector(2,-3){3.1}} 
\put(8.7,-0.2){\vector(-2,-3){3.1}}
\put(4.6,-5.4){\vector(-3,-2){3.2}} 
\put(5.5,-5.4){\vector(3,-2){3.2}}
\put(1.7,-7.7){\vector(1,0){6.4}} 
\put(8.1,-7.9){\vector(-1,0){6.4}}
\put(4.4,-2.8){\vector(-2,-3){3.1}} 
\put(5.7,-2.8){\vector(2,-3){3.1}}
\put(3,-1.1){\ACA?} \put(3.3,-1.35){\f{ats->HB}}
\put(6.2,-1.1){\RCA} \put(6.3,-1.35){\f{leb->HB}}
\put(2.1,-2.8){\ACA} \put(2.3,-2.5){\f{ats->comp}}
\put(7.1,-2.8){\ACA} \put(7.3,-2.5){\f{leb->comp}}
\put(2.1,-4.9){\RCA} \put(2.4,-5.15){\f{HB->ats}}
\put(7.2,-4.9){\RCA} \put(7.35,-5.15){\f{HB->leb}}
\put(3.1,-6.7){\WKL} \put(3.4,-6.4){\f{comp->ats}}
\put(6.3,-6.7){\WKL} \put(6.3,-6.4){\f{comp&leb}}
\put(4.1,-7.6){\ACA? \f{ats->leb}}
\put(4.2,-8.3){\RCA\ \f{leb->ats}}
\end{picture}
\end{center}

The plan of the paper is as follows. The last part of this 
introductory section gives a brief presentation (for a more 
detailed presentation see \cite{sosoa}) of the subsystems we will 
deal with. Section \ref{basic} reviews the basic techniques and 
results about the coding of complete separable metric spaces in 
second order arithmetic. In section \ref{compact} we deal with 
various notions of compact space and introduce Lebesgue and Atsuji 
spaces. In sections \ref{comp-ats} and \ref{comp-leb} we study the 
relationships between the various notions of compactness on one 
side and Atsuji and Lebesgue spaces on the other side: we prove 
the results expressed by the diagonal arrows of our picture (the 
revelsals of the arrows pointing towards ``compact'' are proved in 
\WKL\ and we show that this is necessary). In section \ref{leb&ats} 
we consider the equivalence of Lebesgue and Atsuji spaces and we 
prove the results expressed by the horizontal arrows: the ``easy'' 
direction (i.e.\ Lebesgue implies Atsuji) can be obtained in \RCA, 
while the opposite direction is provable in \ACA\ (we do not know 
whether it is equivalent to \ACA). Section \ref{dist} considers a 
problem arisen during the investigations of section \ref{leb&ats}, 
namely showing that the function distance from a closed set is 
continuous: this requires the even stronger subsystem \PCA.

The formal systems we will consider are, in order of increasing 
strength, \RCA, \WKL, \ACA\ and \PCA. These are all theories which 
use classical first order logic and the language of second order 
arithmetic, which consists of number variables $m,n,\dots$, set 
variables $X,Y,\dots$, primitives $+$, $\cdot$, $0$, $1$, $=$, $<$ 
and $\in$, logical connectives and quantifiers on both sorts of 
variables. Formulas of this language are classified according to 
the number of alternating quantifiers: \S01 formulas have one 
existential number quantifier in front of a matrix containing only 
bounded number quantifiers; arithmetical formulas contain no set 
quantifiers; \S11 formulas have one existential set quantifier in 
front of an arithmetical matrix; \P in formulas are negations of 
\S in formulas.

All systems share a set of basic arithmetical axioms, an induction 
axiom
$$ 0 \in X \land \forall n (n \in X \imply n+1 \in X) \imply 
\forall n (n \in X)$$
and differ by the formulas $\varphi$ allowed in the comprehension 
scheme
$$\exists X\; \forall n (n \in X \sse \varphi(n))$$
or by the presence of other additional axioms.

\RCA\ has comprehension only for \D01 formulas, i.e.\ formulas 
which are provably equivalent both to a \S01 and to a \P01 formula 
(and, for technical reasons, it has also an induction scheme for 
\S01 formulas): this is the base theory for most reverse 
mathematics investigations. \WKL\ extends \RCA\ by adding to it 
weak K\"onig's lemma (i.e.\ K\"onig's lemma for trees consisting 
of sequences of 0's and 1's): this allows for a good theory of 
compactness and continuity. \ACA\ has comprehension for arbitrary 
arithmetical formulas and allows for a good theory of sequential 
convergence. \PCA\ is the strongest system that turns out to be 
needed to prove theorems of ordinary mathematics: \P11 formulas 
are allowed in the comprehension scheme.

A typical reverse mathematics result is the statement that, within 
a weaker base theory (typically \RCA, but see theorems 
\ref{WKLnecHB}, \ref{WKLnecats} and \ref{WKLnecleb}), one of these 
systems is equivalent to some theorem of ordinary mathematics.

In the next sections whenever we begin a definition, lemma or 
theorem by the name of one of these subsystems between parenthesis 
we mean that the definition is given, or the statement provable, 
within that subsystem.

\section{Coding complete separable metric spaces}\label{basic}
\begin{definition}[\RCA]
A (code for a) {\em complete separable metric space $\A$\/} is a 
set $A \subseteq \N$ together with a function $d: A \times A \to 
\R$ such that for all $a,b,c \in A$ we have $d(a,a) = 0$, $d(a,b) 
= d(b,a) \geq 0$ and $d(a,b) \leq d(a,c) + d(c,b)$.

A (code for a) {\em point of $\A$\/} is a sequence $\<a_n: n \in 
\N\>$ of elements of $A$ such that for every $n$ we have 
$d(a_n,a_{n+1}) < 2^{-n}$.
\end{definition}

Within \RCA\ (or any other subsystem of second order arithmetic) 
$\A$ does not formally exist as a set: notations as $x \in \A$ are 
just abbreviations for ``$x$ is a point of $\A$''. Similar 
considerations can be made for the various notions of subsets of 
$\A$ we will introduce.

The metric $d$ can be extended to $\A \times \A$ in an obvious 
way: this extension will still be denoted by $d$ (or by $d_{\A}$ 
when there is danger of confusion) and represents the metric of 
the complete separable metric space. If $x,y \in \A$ are such that 
$d(x,y) = 0$ we identify them and write $x=y$.

A standard example of a complete separable metric space is 
obtained by taking $A = \Q$ with $d$ the usual metric and denoting 
$\widehat{\Q}$ by $\R$. By restricting the code to the rationals 
between $0$ and $1$ we get a code for the closed interval $[0,1]$.

Other important complete separable metric spaces are the Cantor 
space $\Can$ and the Baire space $\Bai$ of infinite sequences 
respectively of $0$'s and $1$'s and of natural numbers with the 
product topology originated by the discrete topologies 
respectively on $\{0,1\}$ and $\N$. Complete metrics compatible 
with these topologies are obtained by setting, whenever $x \neq 
y$, $d(x,y) = 2^{-n}$ where $n$ is least such that $x(n) \neq 
y(n)$. Explicit codings of $\Can$ and $\Bai$ within \RCA\ are 
provided in \cite{sosoa} and \cite{br-thesis}.

\begin{definition}[\RCA]
For every $x \in \A$ and $q \in \R^+$ let $B(x,q)$ denote the {\em 
open ball\/} of center $x$ and radius $q$ in $\A$. This means that 
for every $y \in \A$ we have that $y \in B(x,q)$ if and only if 
$d(x,y) < q$.

A (code for an) {\em open set\/} in $\A$ is a sequence $U = 
\<(a_n,q_n): n \in \N\>$ of elements of $A \times \Q^+$. The 
meaning of this coding is that $U = \bigcup_{n \in \N} B(a_n,q_n)$ 
and hence $x \in U$ if and only if $\exists n\; d(x,a_n) < q_n$. A 
{\em closed set\/} in $\A$ is the complement of an open set, and 
thus is represented by the same code.
\end{definition}

A basic fact about open sets in complete separable metric spaces 
is the following lemma proved in \cite{sosoa}.

\begin{lemma}[\RCA]\label{s01open}
Let $\varphi(x)$ be a \S01 formula such that $x,y \in \A$ and 
$x=y$ imply $\varphi(x) \sse \varphi(y)$. Then there exists an 
open set $U$ in $\A$ such that $x \in U$ if and only if 
$\varphi(x)$ holds.
\end{lemma}

\begin{definition}[\RCA]
For every $x \in \A$ and $q \in \R^+$ let
$$P(x,q) =\set{y \in \A}{0 < d(x,y) < q} = B(x,q) \setminus \{x\}$$
be the {\em punctured ball\/} of center $x$ and radius $q$.
\end{definition}

The formula defining $P(x, q)$ is \S01: hence by lemma 
\ref{s01open} within \RCA\ $P(x,q)$ is an open set.

\begin{definition}[\RCA]
A point $x \in \A$ is {\em isolated\/} if for some $q \in \R^+$ we 
have $P(x,q) = \emptyset$. A complete separable metric space is 
{\em perfect\/} if it does not have isolated points.
\end{definition}

We defined closed sets to be complements of open sets; another 
natural definition can be obtained by viewing a closed set as the 
closure of a countable set.

\begin{definition}[\RCA]
A code for a {\em separably closed set\/} in $\A$ is a sequence $C= 
\<x_n: n \in \N\>$ of points of $\A$. The separably closed set is 
then denoted by $\cl C$ and $x \in \cl C$ if and only if $\forall 
q \in \Q^+\; \exists n\; d(x,x_n) < q$.
\end{definition}

The two notions of closed set we introduced are not equivalent 
within \RCA: their relationship has been studied in depth by Brown 
(\cite{br-thesis,br}), who obtained the following results.

\begin{theorem}[\RCA]\label{sep->clo}
The following are equivalent:\begin{enumerate}
\item \ACA.
\item Every separably closed set in a complete separable metric 
space is closed.
\end{enumerate}\end{theorem}

\begin{theorem}[\RCA]\label{clo->sep}
The following are equivalent:\begin{enumerate}
\item \PCA.
\item Every closed set in a complete separable metric space is 
separably closed.
\end{enumerate}\end{theorem}

If $W$ and $Z$ are open, closed or separably closed sets of $\A$ 
we write $W \subseteq Z$ to mean $\forall x (x \in W \imply x \in 
Z)$. $W = Z$ and $W \nsubseteq Z$ have the obvious meanings.

Continuous functions are coded in second order arithmetic in the 
following way (see \cite{bs1,sosoa}).

\begin{definition}[\RCA]
Let $\A$ and $\B$ be two complete separable metric spaces. A (code 
for a) {\em continuous function from $\A$ to $\B$\/} is a set $\Phi 
\subseteq \N \times A \times \Q^+ \times B \times \Q^+$ such that, 
if we denote by $(a,r)\Phi(b,s)$ the formula $\exists n\; 
(n,a,r,b,s) \in \Phi$, the following properties 
hold:\begin{itemize}
\item $(a,r)\Phi(b,s) \land (a,r)\Phi(b',s') \imply d(b,b') < s + 
s'$;
\item $(a,r)\Phi(b,s) \land d(b,b') + s \leq s' \imply 
(a,r)\Phi(b',s')$;
\item $(a,r)\Phi(b,s) \land d(a,a') + r' \leq r \imply 
(a',r')\Phi(b,s)$;
\item $\forall x \in \A\; \forall q \in \Q^+ \exists (a,r,b,s) 
((a,r)\Phi(b,s) \land d(x,a) < r \land s < q)$.
\end{itemize}
In this situation for every $x \in \A$ there exists a unique $y 
\in \B$ such that $d(y,b) < s$ whenever $d(x,a) < r$ and 
$(a,r)\Phi(b,s)$. This $y$ is denoted by $f(x)$ and is the image 
of $x$ under the function $f$ coded by $\Phi$.
\end{definition}

Sometimes we will need to use continuous functions which are 
defined only on a subset of $\A$. These can be coded omitting the 
last clause in the above definition: their domain consists 
precisely of those $x \in \A$ for which
$$\forall q \in \Q^+ \exists (a,r,b,s) ((a,r)\Phi(b,s) \land 
d(x,a) < r \land s < q).$$

\RCA\ proves (see \cite{sosoa}) that $d$ is a continuous function, 
that the class of continuous function contains the constant 
functions and is closed under the basic arithmetical operations, 
$\max$, $\min$ and composition and that the preimage by a 
continuous function of an open set in $\B$ is an open set in $\A$. 
We will use these facts without explicit mention.

In section \ref{dist} we will show that not all continuous 
functions between complete separable metric spaces which are 
commonly used exist within \RCA: indeed constructing for every 
closed set a code for the continuous function that associates to 
every point its distance from the closed set requires \PCA, while 
the same construction for separably closed sets requires \ACA.

The following versions of Urysohn's lemma and Tietze extension 
theorem are proved in \cite{br-thesis} and \cite{sosoa}.

\begin{theorem}[\RCA]\label{ury}
If $C_0$ and $C_1$ are closed sets in a complete separable metric 
space $\A$ and $C_0 \cap C_1 = \emptyset$ then there exists a 
continuous function $f:\A \to \R$ such that for every $i<2$ and 
$x\in C_i$ we have $f(x)=i$.
\end{theorem}

\begin{theorem}[\RCA]\label{tietze}
If $C$ is a closed set in a complete separable metric space $\A$ 
and $f: C \to \R$ is a continuous function there exists a 
continuous function $g: \A \to \R$ such that $g \restriction C = 
f$, i.e.\ $g(x) = f(x)$ for every $x \in C$.
\end{theorem}

We will consider also uniformly continuous functions. In the 
context of subsystems of second order arithmetic sometimes (e.g.\ 
in \cite{sosoa}) functions which admit a modulus of uniform 
continuity have been most useful. Here we consider the usual (and 
weaker, from the point of view of subsystems of second order 
arithmetic) notion of uniformly continuous function because we 
want to be as close as possible to standard mathematical practice.

\begin{definition}[\RCA]
A continuous function $f: \A \to \B$ is {\em uniformly 
continuous\/} if
$$\forall \varepsilon \in \R^+\; \exists \delta \in \R^+\; \forall 
x,y \in \A (d_{\A}(x,y) < \delta \imply d_{\B}(f(x),f(y)) < 
\varepsilon).$$
\end{definition}

\section{Compact, Lebesgue and Atsuji spaces}\label{compact}
Another important concept is that of compact space. As is the case 
for the notions of closed sets also in this case the usual 
equivalent notions of compact spaces are not equivalent in weak 
subsystems of second order arithmetic. We will consider the 
following notions.

\begin{definition}[\RCA]
A complete separable metric space $\A$ is {\em compact\/} if there 
exists an infinite sequence of finite sequences of points of $\A$ 
$\<\<x_{n,m}: m < i_n\> : n \in \N\>$ such that
$$\forall x \in \A\; \forall n \in \N\; \exists m < i_n\; 
d(x,x_{n,m}) < 2^{-n}.$$
\end{definition}

Notice that our definition of compact space requires more than the 
existence for every $\varepsilon \in \R^+$ of a finite set $B 
\subseteq \A$ such that $\forall x \in \A\; \exists y \in B\; 
d(x,y) < \varepsilon$. \RCA\ does not prove that the latter 
condition implies compactness but \ACA\ does and indeed proves the 
following equivalence.

\begin{lemma}[\ACA]\label{nets}
A complete separable metric space $\A$ is compact if and only if 
for every $\varepsilon \in \R^+$ there exists a finite set $B 
\subseteq A$ such that $\forall a \in A\; \exists b \in B\; d(a,b) 
\leq \varepsilon$.
\end{lemma}
\begin{pf}
Suppose $\<\<x_{n,m}: m < i_n\> : n \in \N\>$ witnesses the 
compactness of $\A$. For any $\varepsilon \in \R^+$ let $n$ be 
such that $2^{-n+1} \leq \varepsilon$. For every $m < i_n$ let 
$b_m \in A$ be such that $d(b_m, x_{n,m}) < 2^{-n}$. Then $B = 
\set{b_m}{m<i_n}$ satisfies $\forall x \in \A\; \exists b \in B\; 
d(x,b) \leq \varepsilon$ and a fortiori the desired property.

For the other direction of the equivalence recall that, within 
\RCA, every finite set can be coded as a natural number. Hence if 
for every $n \in \N$ there exists a finite set $B \subseteq A$ 
such that $\forall a \in A\; \exists b \in B\; d(a,b) \leq 
2^{-(n+1)}$, within \ACA, there exists a function $f$ that to each 
$n$ associates the least finite set $B$ satisfying this 
arithmetical condition. If $\<\<x_{n,m}: m < i_n\> : n \in \N\>$ is 
the sequence of the sequences of the elements of the various 
$f(n)$'s it is immediate to check that it witnesses the 
compactness of $\A$.
\end{pf}

\begin{definition}[\RCA]
A sequence $\U = \<U_n:n \in \N\>$ of open sets in $\A$ is an {\em 
open covering} if for every $x \in \A$ there exists $n$ such that 
$x \in U_n$.

A complete separable metric space $\A$ is {\em Heine-Borel 
compact\/} if for every open covering $\U$ of $\A$ there exists a 
finite covering $\U' \subseteq \U$.
\end{definition}

The results summarized in the next theorem are contained in 
\cite{br-thesis} and \cite{sosoa}.

\begin{theorem}[\RCA]\label{comp->HB}
The following are equivalent:\begin{enumerate}
\item \WKL.
\item Every compact complete separable metric space is Heine-Borel 
compact.
\item The closed interval $[0,1]$ is Heine-Borel compact.
\item If $\A$ is a compact complete separable metric space and 
$f:\A \to \R$ a continuous function then $f$ attains a minimum.
\end{enumerate}\end{theorem}

We also have the following result which is a corollary of the 
proof of one of the implications of the above theorem.

\begin{lemma}[\RCA]\label{HB->min}
If $\A$ is a Heine-Borel compact complete separable metric space 
and $f:\A \to \R$ a continuous function then $f$ attains a 
minimum.
\end{lemma}
\begin{pf}
It suffices to inspect the proof of $(1) \implies (4)$ of the 
previous theorem in \cite{sosoa} and notice that \WKL\ is only used 
to deduce Heine-Borel compactness from compactness.
\end{pf}

The following equivalence is essentially due to Brown 
(\cite{br-thesis}), but we prove a slightly different result that 
will be useful in the proof of theorem \ref{ats->comp}.

\begin{theorem}[\WKL]\label{HB->comp}
The following are equivalent:\begin{enumerate}
\item \ACA.
\item Every Heine-Borel compact complete separable metric space is 
compact.
\item Every perfect Heine-Borel compact complete separable metric 
space is compact.
\end{enumerate}\end{theorem}
\begin{pf}
The equivalence between (1) and (2) is proved in \cite{br-thesis} 
within \RCA. Since (2) implies (3) is obvious it suffices to prove 
that (3) implies (1): to this end we modify slightly Brown's proof 
of (2) implies (1).

It is well-known that \ACA\ is equivalent over \RCA\ (and, a 
fortiori, over \WKL) to the statement that the range of every 
one-to-one function from $\N$ to $\N$ exists. Fix $f: \N \to \N$ 
one-to-one. We want to define a code for a complete separable 
metric space $\A$ homeomorphic to
$$\{(0,0)\} \cup \bigcup_{k \in \N} \left(\{2^{-f(k)}\} \times 
[0,2^{-f(k)}]\right) \subseteq \R^2.$$
Since the range of $f$ is not available we cannot use 
$\set{(2^{-f(k)},q)}{0 \leq q \leq 2^{-f(k)}}$ as a code. This 
problem can be overcome by defining
$$A = \set{(k,q) \in \N \times \Q^+}{q \leq 2^{-f(k)}}$$
and letting
$$d((k,q),(k',q')) = \max(|2^{-f(k)} - 2^{-f(k')}|, |q - q'|)$$
(we are using a metric different from, but equivalent to, the 
usual metric on $\R^2$).

$\A$ is clearly perfect and we claim that it is also Heine-Borel 
compact. To see this suppose $\U = \<U_n : n \in \N\>$ is an open 
covering of $\A$. For some $n_0$ we have that $(0,0) \in U_{n_0}$ 
and hence there exists $m$ such that if $f(k) > m$ then $(k,y) \in 
U_{n_0}$ for every $y \in [0,2^{-f(k)}]$. There are only finitely 
many $k$'s such that $f(k) \leq m$ and we can define $n_1 = \max 
(\{n_0\} \cup \{g(k) : f(k) \leq m\})$, where $g:\N \to \N$ is such 
that for every $k$ we have that $\<U_n : n < g(k)\>$ is a covering 
of $\{2^{-f(k)}\} \times [0,2^{-f(k)}]$ ($g$ exists within \WKL\ by 
the uniform version of the Heine-Borel compactness of $[0,1]$ 
proved in \cite{sosoa}). We have that $\<U_n : n < n_1\>$ is a 
finite subcovering of $\U$.

Therefore (3) implies that $\A$ is compact: let $\<\<x_{n,m}: m < 
i_n\> : n \in \N\>$ be such that
$$\forall x \in \A\; \forall n \in \N\; \exists m < i_n\; 
d(x,x_{n,m}) < 2^{-n}.$$
Every $x_{n,m}$ represents some $(2^{-f(k_{n,m})}, y_{n,m})$, 
whose actual code is $(k_{n,m}, y_{n,m})$. It is easy to check 
that for every $n$ we have
$$\exists k\; f(k) = n \text{\quad iff \quad} \exists m < 
i_{n+1}\; f(k_{n+1,m}) = n.$$
By recursive comprehension the range of $f$ exists.
\end{pf}

The reader may notice that the preceding theorem has been proved 
within \WKL, while most reverse mathematics results are proved 
within \RCA. The use of a stronger base theory is indeed necessary 
to prove that statement (3) implies \ACA\ in this theorem, as we 
are now going to show, and the same situation will occur also for 
other results we will obtain later.

We need to formalize within \RCA\ the fact that every perfect 
complete separable metric space has a closed subset which is 
homeomorphic to $\Can$ (with the metric described in section 
\ref{basic}). This amounts essentially to Exercise 3D.15 in 
\cite{mos} (which does not mention \RCA).

\begin{theorem}[\RCA]\label{embeds}
For every perfect complete separable metric space $\A$ there 
exists a sequence $\<B(a_s,q_s) :s \in \Seq\>$ of open balls with 
$a_s \in A$ and $q_s \in \R^+$ such that:\begin{enumerate}
\item $\forall s \in \Seq\; \forall i<2\; d(a_s, a_{s \conc 
\<i\>}) + q_{s \conc \<i\>}< q_s$;
\item $\forall s \in \Seq\; d(a_{s\conc \<0\>}, a_{s \conc \<1\>}) 
< q_{s \conc \<0\>} + q_{s \conc \<1\>}$;
\item $\forall s \in \Seq\; q_s \leq 2^{-\lh(s)}$.
\end{enumerate}

Using $\<B(a_s,q_s) :s \in \Seq\>$ we can define an embedding 
(i.e.\ a continuous function which is injective and has 
continuous inverse) $\varphi : \Can \to \A$ such that the range of 
$\varphi$ is closed in $\A$. Moreover $\varphi$ is Lipschitz with 
constant $1$, i.e.\ $d_{\A}(\varphi(x), \varphi(y)) \leq 
d_{\Can}(x,y)$ for every $x,y \in \Can$.
\end{theorem}
\begin{pf}
Define $a_s$ and $q_s$ by recursion on $\lh(s)$. Let $a_{\<\>}$ be 
an element of $A$ and $q_{\<\>} = 1$. Assuming we have defined 
$a_s$ and $q_s$, since $\A$ is perfect $B(a_s,q_s)$ contains at 
least two distinct elements $a_{s \conc \<0\>}$ and $a_{s \conc 
\<1\>}$. For every $i<2$ let
$$q_{s \conc \<i\>} = \min \left\{\tfrac 12 (q_s - d(a_s, a_{s 
\conc \<i\>})), \tfrac 13 d(a_{s \conc \<0\>}, a_{s \conc 
\<1\>})\right\}.$$

$\varphi$ is defined by letting, for every $x \in \Can$, 
$\varphi(x) = \<a_{x \restriction n} : n \in \N\>$, where $x \restriction n$ 
is the initial segment of $x$ of length $n$. The properties of 
$\<B(a_s,q_s) :s \in \Seq\>$ imply that $\varphi(x)$ is a point of 
$\A$ and that $\varphi$ is injective and Lipschitz. The complement 
of the range of $\varphi$ is $U = \set{x \in \A}{\exists n\; 
\forall s \in 2^n\; d(x,a_s) > q_s}$ which is open by lemma 
\ref{s01open}. We leave to the reader the routine details (which 
involve the details of the coding of $\Can$) of the definition of 
codes for $\varphi$ and its inverse.
\end{pf}

\begin{lemma}[\RCA]\label{perfHB}
The following are equivalent:\begin{enumerate}
\item \WKL.
\item There exists a complete separable metric space which is 
perfect and Heine-Borel compact.
\end{enumerate}\end{lemma}
\begin{pf}
(1) implies (2) follows immediately from theorem \ref{comp->HB} 
and the fact that \RCA\ proves that $[0,1]$ is perfect.

To prove that (2) implies (1) let $\A$ be a perfect Heine-Borel 
compact complete separable metric space. Let $\<B(a_s,q_s) :s \in 
\Seq\>$ and $\varphi$ be given by theorem \ref{embeds} and denote 
by $U$ be the complement of the range of $\varphi$.

Now let $T \subseteq \Seq$ be a binary tree with no paths. 
Consider the collection of open sets $\U = \{U\} \cup 
\set{B(a_s,q_s)}{s \notin T}$. $\U$ is a covering of $\A$ because 
$T$ has no paths. Let $\U' \subseteq \U$ be a finite subcovering: 
only for finitely many $s$ we have $B(a_s,q_s) \in \U'$. Since 
$\U'$ is a covering every $t \in T$ has an extension $s$ such that 
$B(a_s,q_s) \in \U'$ and for each $s$ there are only finitely many 
such $t$, this entails that $T$ is finite.
\end{pf}

\begin{theorem}\label{WKLnecHB}
Statement (3) of theorem \ref{HB->comp} does not imply \ACA\ in 
any theory stronger than \RCA\ and properly weaker than \WKL.
\end{theorem}
\begin{pf}
Let $\system{T}$ be a theory stronger than \RCA\ and properly 
weaker than \WKL\ and let $\frak M$ be a model of $\system{T}$ in 
which \WKL\ fails. By lemma \ref{perfHB} in $\frak M$ there are no 
perfect complete separable metric spaces which are Heine-Borel 
compact and hence statement (3) of theorem \ref{HB->comp} is 
vacuously true. Since $\frak M$ is not a model of \ACA\ 
$\system{T}$ does not prove that (3) implies \ACA.
\end{pf}

Our main goal is to explore the relationships among the following 
notions and between them and the notions of compactness we just 
introduced.

\begin{definition}[\RCA]
A complete separable metric space $\A$ is {\em Atsuji} if every 
continuous function $f: \A \to \B$ (where $\B$ is an arbitrary 
complete separable metric space) is uniformly continuous.
\end{definition}

\begin{definition}[\RCA]
A complete separable metric space $\A$ is {\em Lebesgue} if for 
every open covering $\U = \<U_n:n \in \N\>$ of $\A$ there exists 
$q \in \R^+$ such that
$$\forall x \in \A\; \exists n \in \N\; B(x,q) \subseteq U_n.$$
$q$ is called a {\em Lebesgue number\/} for $\U$.
\end{definition}

\begin{remark}\label{example}
The set of natural numbers $\N$ with the usual metric is a 
complete separable metric space which is Atsuji (for every 
$\varepsilon > 0$, $\delta = 1$ suffices in the definition of 
uniform continuity) and Lebesgue ($1$ is a Lebesgue number for 
every covering of $\N$) but not compact.

Another example of a Lebesgue and Atsuji non compact complete 
separable metric space is obtained by taking $\set{e_n}{n \in \N}$ 
to be an orthonormal basis for an infinite dimensional separable 
real Hilbert space and considering $\{0\} \cup \set{2^{-m}e_n}{m,n 
\in \N}$.
\end{remark}

\section{Compact and Atsuji}\label{comp-ats}
One direction of the relationship between compact spaces and 
uniform continuity has been already explored by Brown and Simpson. 
The following is the statement in our terminology of the main 
results they obtained.

\begin{theorem}[\RCA]\label{comp->ats}
The following are equivalent:\begin{enumerate}
\item \WKL.
\item Every complete separable metric space which is compact is 
Atsuji.
\item The closed interval $[0,1]$ is Atsuji.
\end{enumerate}\end{theorem}
\begin{pf}
(1) implies (2) is proved in \cite{br-thesis} and \cite{sosoa}. (2) 
implies (3) holds because in \RCA\ it is easy to show that $[0,1]$ 
is compact. (3) implies (1) is proved in \cite{sosoa}.
\end{pf}

We also have the following result which is a corollary of the 
proof of one of the implications of the above theorem.

\begin{lemma}[\RCA]\label{HB->ats}
Every complete separable metric space which is Heine-Borel compact 
is Atsuji.
\end{lemma}
\begin{pf}
It suffices to inspect the proof of implication $(1) \implies (2)$ 
of the previous theorem and notice that \WKL\ is only used to 
deduce Heine-Borel compactness from compactness.
\end{pf}

Remark \ref{example} shows that not all Atsuji spaces are compact. 
However a perfect space which is Atsuji is compact. \ACA\ is 
needed to prove this result.

\begin{theorem}[\WKL]\label{ats->comp}
The following are equivalent:\begin{enumerate}
\item \ACA.
\item Every complete separable metric space which is perfect and 
Atsuji is compact.
\end{enumerate}\end{theorem}
\begin{pf}
$(1) \implies (2)$. We reason in \ACA\ and suppose that $\A$ is a 
perfect complete separable metric space which is not compact. We 
will show that $\A$ is not Atsuji.

Since $\A$ is not compact by lemma \ref{nets} there exists 
$\varepsilon \in \R^+$ such that for no finite $B \subseteq A$ we 
have that for all $a \in A$ there exists $b \in B$ such that 
$d(a,b) \leq \varepsilon$. Using this fact we can define by 
recursion a sequence $\<a^0_n: n \in \N\>$ of elements of $A$ such 
that $n \neq m \imply d(a^0_n,a^0_m) > \varepsilon$. Since $\A$ is 
perfect for every $n$ there exists $a^1_n \in A$ such that $a^1_n 
\in P(a^0_n,2^{-n-1}\varepsilon)$. Using the triangle inequality we 
get that $n \neq m \imply d(a^1_n,a^1_m) > \frac\varepsilon4$. 
Setting $C_i = \<a^i_n: n \in \N\>$ (for $i= 0,1$) these facts 
entail that $C_0 \cap C_1 = \emptyset$ and $\cl {C_i} = C_i$. In 
other words, $C_0$ and $C_1$ code two disjoint separably closed 
sets. By theorem \ref{sep->clo} each $C_i$ is closed and we can 
apply theorem \ref{ury} to get a continuous function $f : \A \to 
\R$ such that $f(C_i) = \{i\}$. To see that $f$ is not uniformly 
continuous fix $\delta \in \R^+$: if $2^{-n-1}\varepsilon \leq 
\delta$ we have $d(a^0_n,a^1_n) < \delta$ but $|f(a^0_n) - 
f(a^1_n)| = 1$.

$(2) \implies (1)$. We will use theorem \ref{HB->comp}: it suffices 
to prove that if $\A$ is Heine-Borel compact and perfect then $\A$ 
is compact. This follows immediately from lemma \ref{HB->ats} and 
(2).
\end{pf}

\begin{remark}\label{ats->HB}
Combining theorems \ref{ats->comp} and \ref{comp->HB} we have a 
proof within \ACA\ that every complete separable metric space 
which is perfect and Atsuji is Heine-Borel compact. We do not know 
whether \ACA\ is necessary to prove this statement.
\end{remark}

We will prove that \WKL\ is necessary to obtain the equivalence of 
theorem \ref{ats->comp} by the same argument we used to prove 
theorem \ref{WKLnecHB}.

\begin{lemma}[\RCA]\label{perfats}
The following are equivalent:\begin{enumerate}
\item \WKL.
\item There exists a complete separable metric space which is 
perfect and Atsuji.
\end{enumerate}\end{lemma}
\begin{pf}
(1) implies (2) follows immediately from theorem \ref{comp->ats}.

To prove that (2) implies (1) let $\A$ be a perfect Atsuji 
complete separable metric space and $\varphi$ be given by theorem 
\ref{embeds}. Denote by $C$ and $\psi$ respectively the range and 
the inverse of $\varphi$. We will show that every continuous 
function $f: \Can \to \R$ is uniformly continuous. This implies 
\WKL\ by a simplified version of the argument used in \cite{sosoa} 
to prove $(3) \implies (1)$ of theorem \ref{comp->ats} (that 
argument uses only functions from $[0,1]$ to $\R$ and the explicit 
embedding of $\Can$ into $[0,1]$ given by Cantor middle-third 
set).

Let $f: \Can \to \R$ be continuous and define $g: C \to \R$ by 
setting $g = f \circ \psi$. By theorem \ref{tietze} let $h: \A \to 
\R$ be continuous such that $h \restriction C = g$. Since $\A$ is 
Atsuji $h$ is uniformly continuous. To show that $f$ is uniformly 
continuous fix $\varepsilon > 0$ and let $\delta > 0$ be such that 
for all $y,y' \in \A$ if $d_{\A}(y,y') < \delta$ then $|h(y) - 
h(y')| < \varepsilon$. If $x,x' \in \Can$ are such that 
$d_{\Can}(x,x') < \delta$ then, since $\varphi$ is Lipschitz with 
constant $1$, $d_{\A}(\varphi(x), \varphi(x')) < \delta$ and hence 
$|h(\varphi(x)) - h(\varphi(x'))| < \varepsilon$. But 
$h(\varphi(x)) = f(x)$ and $h(\varphi(x')) = f(x')$, so that the 
uniform continuity of $f$ is established.
\end{pf}

\begin{theorem}\label{WKLnecats}
Statement (2) of theorem \ref{ats->comp} does not imply \ACA\ in 
any theory stronger than \RCA\ and properly weaker than \WKL.
\end{theorem}
\begin{pf}
Repeat the argument of the proof of theorem \ref{WKLnecHB} using 
lemma \ref{perfats}.
\end{pf}

\section{Compact and Lebesgue}\label{comp-leb}
We now explore the relationship between compact spaces and 
Lebesgue numbers.

\begin{lemma}[\WKL]\label{comp->leb}
Every complete separable metric space which is compact is 
Lebesgue.
\end{lemma}
\begin{pf}
We reason in \WKL. Let $\A$ be a compact complete separable metric 
space and $\U$ an open covering of $\A$. Each element of $\U$ is 
union of open balls with center in $A$ and rational radius: since 
a Lebesgue number for the covering consisting of these open balls 
is also a Lebesgue number for the original covering we can assume 
that each element of $\U$ is actually such an open ball. By 
theorem \ref{comp->HB} $\A$ is Heine-Borel compact and there exists 
a finite subcovering of $\U$. Since a Lebesgue number for any 
subcovering is a Lebesgue number also for the original covering we 
can assume that $\U$ is finite and has the form $\<B(a_n,r_n): 
n<k\>$ with $a_n \in A$ and $r_n \in \Q^+$ for every $n<k$.

For every $n<k$ let $f_n: \A \to \R$ be the continuous function 
defined by $f_n(x) = \max(0,r_n - d(a_n,x))$. Let $f:\A \to \R$ be 
the continuous function defined by $f(x) = \max 
\set{f_n(x)}{n<k}$. Since $\U$ is a covering for every $x \in \A$ 
there exists $n<k$ such that $f_n(x) > 0$ and hence $f(x) > 0$. By 
\ref{comp->HB} since $\A$ is compact $f$ attains a minimum $q \in 
\R^+$. We claim that $q$ is a Lebesgue number for $\U$.

To prove the claim let $x \in \A$: since $f(x) \geq q$ for some 
$n<k$ we have $r_n - d(a_n,x) = f_n(x) \geq q$ which implies 
$B(x,q) \subseteq B(a_n,r_n)$ completing the proof of the claim 
and of the lemma.
\end{pf}

\begin{remark}
The functions $f_n$ used in the proof of the preceding lemma 
compute a lower bound for the distance from the complement of 
$B(a_n,r_n)$ (the latter is the function used in textbook proofs 
of this result). This suffices to prove that if $f_n(x) \geq q$ 
then $B(x,q) \subseteq B(a_n,r_n)$, which is all is needed to 
complete the argument. In theorem \ref{cont-d} we will show that 
\WKL\ does not suffice to prove that the function computing the 
actual distance from the complement of $B(a_n,r_n)$ exists.
\end{remark}

The following results will be used in the proofs of theorems 
\ref{comp&leb} and \ref{leb->comp} but are also interesting in 
their own right.

\begin{lemma}[\RCA]\label{HB->leb}
Every complete separable metric space which is Heine-Borel compact 
is Lebesgue.
\end{lemma}
\begin{pf}
It suffices to repeat the proof of the above lemma using lemma 
\ref{HB->min}.
\end{pf}

\begin{theorem}[\RCA]\label{leb->HB}
Every complete separable metric space which is perfect and 
Lebesgue is Heine-Borel compact.
\end{theorem}
\begin{pf}
Let $\U$ be an open covering of $\A$. Since every open subset of 
$\A$ is union of open balls with center in $A$ and rational radius 
we may assume that $\U$ has the form $\<B(a_n,r_n): n \in \N\>$ 
with $a_n \in A$ and $r_n \in \Q^+$ for every $n \in \N$.

For every $n \in \N$ and $b \in A$ let $q_{n,b} = \min (r_n - 
d(a_n,b), 2^{-n}) \in \R$. Let $V_{n,b} = P(b,q_{n,b})$. Notice 
that the definition of $q_{n,b}$ implies that $V_{n,b} \subseteq 
B(a_n,r_n)$ and that if $b \notin B(a_n,r_n)$ we have $q_{n,b} \leq 
0$ and hence $V_{n,b} = \emptyset$.

We claim that $\V = \<V_{n,b}: n \in \N, b \in A\>$ is a covering 
of $\A$. To see this let $x \in \A$: since $\U$ is a covering, for 
some $n \in \N$ we have $d(a_n,x) < r_n$ and hence there exists 
$\varepsilon \in \R^+$ such that $d(a_n,x) \leq r_n - 
2\varepsilon$. Since $\A$ is perfect and hence $x$ is not isolated 
there exists $b \in A$ such that $0 < d(b,x) < \min (\varepsilon, 
2^{-n})$: this implies $d(b,a_n) < r_n - \varepsilon$. Since 
$d(b,x) < \min (\varepsilon, 2^{-n}) < \min (r_n - d(a_n,b), 
2^{-n}) = q_{n,b}$ we have $x \in V_{n,b}$. This completes the 
proof of the claim.

Since $\A$ is Lebesgue there exists a Lebesgue number $q \in \R^+$ 
for $\V$. Let $k \in \N$ be such that $2^{-k} < q$. We now prove 
that $\<B(a_n,r_n): n < k\>$ is a finite subcovering of $\U$, 
thereby establishing the lemma.

To see that $\<B(a_n,r_n): n < k\>$ is a covering of $\A$ let $x 
\in \A$: by definition of Lebesgue number there exist $n \in \N$ 
and $b \in A$ such that $B(x,q) \subseteq V_{n,b}$. Thus $b \notin 
B(x,q)$ and $x \in V_{n,b}$ which imply $q \leq d(b,x) < q_{n,b}$. 
Therefore $2^{-k} < q < q_{n,b} \leq 2^{-n}$ which entails $n<k$. 
Since $x \in V_{n,b} \subseteq B(a_n,r_n)$ the proof is complete.
\end{pf}

The following are our reverse mathematics results on the 
relationship between Lebesgue spaces and compactness.

\begin{theorem}[\RCA]\label{comp&leb}
The following are equivalent:\begin{enumerate}
\item \WKL.
\item Every complete separable metric space which is compact is 
Lebesgue.
\item The closed interval $[0,1]$ is Lebesgue.
\end{enumerate}\end{theorem}
\begin{pf}
(1) implies (2) is lemma \ref{comp->leb}. (2) implies (3) is 
immediate. To prove (3) implies (1) use theorem \ref{comp->HB} and 
notice that $[0,1]$ is perfect: it follows from theorem 
\ref{leb->HB} that it is Heine-Borel compact.
\end{pf}

\begin{theorem}[\WKL]\label{leb->comp}
The following are equivalent:\begin{enumerate}
\item \ACA.
\item Every complete separable metric space which is perfect and 
Lebesgue is compact.
\end{enumerate}\end{theorem}
\begin{pf}
(1) implies (2) follows by theorem \ref{leb->HB} and $(1) \implies 
(2)$ of theorem \ref{HB->comp}.

To prove that (2) implies (1) we use theorem \ref{HB->comp}: we 
suppose that $\A$ is perfect and Heine-Borel compact and show, 
using (2), that it is compact. This is immediate using lemma 
\ref{HB->leb}.
\end{pf}

Also in this case we are able to prove that \WKL\ is necessary to 
obtain the equivalence of theorem \ref{leb->comp}.

\begin{lemma}[\RCA]\label{perfleb}
The following are equivalent:\begin{enumerate}
\item \WKL.
\item There exists a complete separable metric space which is 
perfect and Lebesgue.
\end{enumerate}\end{lemma}
\begin{pf}
(1) implies (2) follows immediately from lemma \ref{comp->leb}.

To prove that (2) implies (1) we could give a proof similar to the 
proofs of the corresponding statement in lemmas \ref{perfHB} and 
\ref{perfats}, but this is not necessary: combining theorem 
\ref{leb->HB} and lemma \ref{perfHB} we obtain an immediate proof.
\end{pf}

\begin{theorem}\label{WKLnecleb}
Statement (2) of theorem \ref{leb->comp} does not imply \ACA\ in 
any theory stronger than \RCA\ and properly weaker than \WKL.
\end{theorem}
\begin{pf}
Repeat the argument of the proof of theorem \ref{WKLnecHB} using 
lemma \ref{perfleb}.
\end{pf}

\section{Atsuji and Lebesgue}\label{leb&ats}
In \cite{bee1} and \cite{bee2} Beer remarks that the notions of 
Atsuji space and Lebesgue space are equivalent. The simmetries of 
theorems \ref{comp->ats} and \ref{comp&leb}, which show that both 
notions can be derived from compactness in \WKL, and of theorems 
\ref{ats->comp} and \ref{leb->comp}, which show that both notions 
imply compactness for perfect spaces in \ACA, may suggest that 
this equivalence should be provable in a rather weak subsystems. 
This, as the next theorem shows, is indeed the case for one 
direction of the equivalence.

\begin{theorem}[\RCA]\label{leb->ats}
Every complete separable metric space which is Lebesgue is Atsuji.
\end{theorem}
\begin{pf}
Let $\A$ be a Lebesgue complete separable metric space, $\B$ a 
complete separable metric space and $f: \A \to \B$ a continuous 
function. Fix $\varepsilon \in \R^+$. For every $b \in B$ let $U_b 
= f^{-1}(B(b, \frac\varepsilon2))$: the continuity of $f$ implies 
that $\U = \<U_b : b \in B\>$ is an open covering of $\A$. Let 
$\delta \in \R^+$ be a Lebesgue number for $\U$.

Suppose that $x,y \in \A$ are such that $d_{\A}(x,y) < \delta$: 
this means that $x,y \in B(x,\delta)$. Since $\delta$ is a 
Lebesgue number for $\U$ there exists $b \in B$ such that 
$B(x,\delta) \subseteq U_b$. Therefore $f(x), f(y) \in B(b, 
\frac\varepsilon2)$ and hence $d_{\B}(f(x),f(y)) < \varepsilon$. 
This completes the proof of the uniform continuity of $f$.
\end{pf}

The reverse implication appears to be harder to prove and we 
present a proof of it within \ACA. We do not know whether it is 
provable in a weaker system. Another clue of the difficulties 
involved in proving this implication (and an earlier asymmetry 
between Atsuji and Lebesgue) derives from the fact that we are 
unable to prove the analogue of theorem \ref{leb->HB} with Atsuji 
in place of Lebesgue (see remark \ref{ats->HB}).

The basic tool for our proof that Atsuji spaces are Lebesgue is 
the notion of $\varepsilon$-witness.

\begin{definition}[\RCA]
Let $\U = \<U_n : n \in \N\>$ be an open covering of the complete 
separable metric space $\A$ and $\varepsilon \in \R^+$. We say 
that $x \in \A$ is an {\em $\varepsilon$-witness for $\U$} if for 
every $n$ we have $B(x,\varepsilon) \nsubseteq U_n$, i.e.\ if for 
every $n$ there exists $y \in \A$ such that $y \in 
B(x,\varepsilon)$ and $y \notin U_n$.
\end{definition}

\begin{lemma}[\RCA]\label{wit}
Let $\U = \<U_n : n \in \N\>$ be an open covering of the complete 
separable metric space $\A$. The following properties are 
equivalent:\begin{enumerate}
\item $\U$ has no Lebesgue number.
\item For every $\varepsilon \in \R^+$ there exists $x \in \A$ 
which is an $\varepsilon$-witness for $\U$.
\item For every $\varepsilon \in \R^+$ there exists $a \in A$ 
which is an $\varepsilon$-witness for $\U$.
\end{enumerate}\end{lemma}
\begin{pf}
The equivalence between (1) and (2) follows immediately from the 
definitions. (3) implies (2) is trivial.

To prove that (2) implies (3) suppose that (2) holds and let 
$\varepsilon \in \R^+$. Let $x \in \A$ be an 
$\frac\varepsilon2$-witness for $\U$ and let $a \in A$ be such 
that $d(x,a) < \frac\varepsilon2$. Then $B(a,\varepsilon) 
\supseteq B(x,\frac\varepsilon2)$ and hence $a$ is an 
$\varepsilon$-witness for $\U$.
\end{pf}

\begin{remark}
The formula asserting that $x$ is an $\varepsilon$-witness for 
$\U$ is of the form $\forall n\; \exists y\; 
\psi(n,y,\varepsilon,x)$ with $\psi$ arithmetical, i.e.\ it is an 
essentially \S11 formula (which in \ACA\ is not even provably 
equivalent to a \S11 formula).

Actually it is easy to see that the set of $\varepsilon$-witnesses 
is a $G_\delta$ (countable intersection of open sets) in $\A$ and 
hence is definable by a \P02 formula. The obvious way of doing 
this requires the set $\set{x}{d(x,\A \setminus U_n) < 
\varepsilon}$ to be open, which is a consequence of the continuity 
of the map $x \mapsto d(x,\A \setminus U_n)$: in theorem 
\ref{cont-d} we will show that these two statements are equivalent 
to \PCA\ and hence not available in \ACA.
\end{remark}

In view of the preceding remark the notion of 
$\varepsilon$-witness appears inadequate for a proof in \ACA: we 
need to modify it by using an arithmetical definition, much more 
manageable within \ACA.

\begin{definition}[\RCA]
Let $\U = \<U_n : n \in \N\>$ be an open covering of the complete 
separable metric space $\A$ and $\varepsilon \in \R^+$. We say 
that $x \in \A$ is a {\em strong $\varepsilon$-witness for $\U$} if 
for every $n$ there exists $b \in A$ such that $b \in 
B(x,\varepsilon)$ and $b \notin U_n$.
\end{definition}

\begin{remark}
Notice that not every $\varepsilon$-witness for $\U$ is a strong 
$\varepsilon$-witness. Indeed there exist complete separable 
metric spaces $\A$ and open coverings $\U$ of $\A$ such that for 
every $\varepsilon \in \R^+$ small enough there exist 
$\varepsilon$-witnesses but no strong $\varepsilon$-witnesses for 
$\U$. An example consists, for $\alpha \in \R^+ \setminus \Q$ with 
$\alpha < 1/2$, of the space $\A = \bigcup_{n\in \N} [n-\alpha 
2^{-n}, n+\alpha 2^{-n}] \subset \R$ coded by $A = \A \cap \Q$ with 
the covering $\U = \set{(n-\alpha 2^{-n}, n+\alpha 2^{-n}], 
[n-\alpha 2^{-n}, n+\alpha 2^{-n})}{n \in \N}$.

This shows that lack of Lebesgue number does not imply existence 
of strong $\varepsilon$-witnesses and we cannot replace strong 
$\varepsilon$-witness in place of $\varepsilon$-witness in the 
statement of lemma \ref{wit}. Nevertheless in the proof of the next 
theorem we will change the open covering we deal with so that we 
can use strong $\varepsilon$-witnesses.
\end{remark}

\begin{theorem}[\ACA]\label{ats->leb}
Every Atsuji complete separable metric space is Lebesgue.
\end{theorem}
\begin{pf}
Suppose $\A$ is a complete separable metric space which is not 
Lebesgue and let $\U$ be an open covering of $\A$ which has no 
Lebesgue number. Starting from $\U$ we will construct a continuous 
function $f: \A \to \R$ which is not uniformly continuous, thereby 
showing that $\A$ is not Atsuji.

The first step in our construction is to replace $\U$ by a finer 
open covering $\V$ which not only has no Lebesgue number but for 
every $\varepsilon \in \R^+$ has a strong $\varepsilon$-witness. 
First of all we may assume that $\U$ consists of open balls with 
center in $A$ and rational radius: let $\U = \<B(a_n,r_n): n \in 
\N\>$. Define $\V = \set{B(a_n,s)}{n\in \N \land  s \in \Q^+ \land 
s < r_n}$. It is straightforward to check that $\V$ is an open 
covering of $\A$.

We claim that for every $\varepsilon \in \R^+$ there exists $a \in 
A$ which is a strong $\varepsilon$-witness for $\V$. To prove the 
claim fix $\varepsilon$ and let, by lemma \ref{wit}, $a \in A$ be 
an $\varepsilon$-witness for $\U$. Now fix $B(a_n,s) \in \V$: 
since $a$ is an $\varepsilon$-witness for $\U$ there exists $y \in 
B(a,\varepsilon) \setminus B(a_n,r_n)$. Let $b \in A$ be such that 
$d(b,y) < \min \{\varepsilon - d(a,y), r_n-s\}$. Then it is 
immediate to check that $b \in B(a,\varepsilon)$ and $b \notin 
B(a_n,s)$, thereby showing that $a$ is a strong 
$\varepsilon$-witness for $\V$ and establishing the claim.

Observe that if $a$ is a (strong) $\varepsilon$-witness then 
$P(a,\varepsilon) \neq \emptyset$ and hence there exists $b \in A$ 
such that $b \in P(a,\varepsilon)$. Therefore, using the fact that 
being a strong $\varepsilon$-witness for $\V$ is an arithmetical 
property, within \ACA\ we can construct a sequence $\<(b^0_m, 
b^1_m): m \in \N\>$ of pairs of elements of $A$ such that for 
every $m$ $b^0_m$ is a strong $2^{-m}$-witness for $\V$ and $b^1_m 
\in P(b^0_m, 2^{-m})$.

The following fact about sequences of $b^i_m$'s will be useful in 
the remainder of the proof.

\begin{sublemma}[\RCA]\label{sub}
Suppose $\<x_k: k \in \N\>$ is a sequence of elements of $A$ such 
that for some unbounded function $g: \N \to \N$ we have that for 
every $k$ $x_k$ is $b^i_{g(k)}$ for some $i < 2$. Then $\lim_{k \to 
\infty} x_k$ does not exist.
\end{sublemma}
\begin{pf}
Suppose $x = \lim_{k \to \infty} x_k$. Since $\V$ is an open 
covering of $\A$ there exist $B(a_n,s) \in \V$ and $\varepsilon 
\in \R^+$ such that $B(x,\varepsilon) \subseteq B(a_n,s)$. Since 
the $x_k$'s converge to $x$ and $g$ is unbounded there exists $k$ 
such that $d(x_k,x) < \frac\varepsilon3$ and $2^{-g(k)} < 
\frac\varepsilon3$. Now it is easy to check that 
$B(b^0_{g(k)},2^{-g(k)}) \subseteq B(x,\varepsilon) \subseteq 
B(a_n,s)$, contradicting the fact that $b^0_{g(k)}$ is a 
$2^{-g(k)}$-witness for $\V$.
\end{pf}

The sublemma implies that for every $i<2$ the sequence $\<b^i_m: m 
\in \N\>$ does not contain infinitely many repetitions of the same 
element of $A$. Hence we can define by recursion a strictly 
increasing function $h: \N \to \N$ by setting $h(0) = 0$ and 
$h(n+1) =$ the least $k$ such that for all $m \geq k$ and all $i,j 
<2$ $b^i_m \neq b^j_{h(n)}$. The definition of $h$ implies that if 
we let $C_i = \set{b^i_{h(n)}}{n \in \N}$ then $C_0 \cap C_1 = 
\emptyset$.

Another consequence of the sublemma is that if a sequence of 
elements of $C_i$ converges it is eventually constant. This means 
$\cl{C_i} = C_i$, i.e. that $C_0$ and $C_1$ are separably closed. 
Exactly as in the proof of $(1) \implies (2)$ in theorem 
\ref{ats->comp} we use theorems \ref{sep->clo} and \ref{ury} to 
construct a continuous function $f: \A \to \R$ such that $f(C_i) = 
\{i\}$. To see that $f$ is not uniformly continuous fix $\delta \in 
\R^+$: let $n$ be such that $2^{-h(n)} \leq \delta$: then 
$d(b^0_{h(n)},b^1_{h(n)}) < \delta$ but $|f(b^1_{h(n)}) - 
f(b^0_{h(n)})| = 1$.
\end{pf}

\section{The continuity of the distance from a closed 
set}\label{dist}
In this section we study the function that computes the distance 
of points of a complete separable metric space from a fixed closed 
set: the definition of this function involves a greatest lower 
bound and it is well-known that the existence of $\inf$'s and 
$\sup$'s is equivalent to \ACA\ (see \cite{sosoa}: this is indeed 
one of the very first reverse mathematics results obtained by 
Friedman). However we show that \ACA\ does not suffice to prove 
the continuity of the distance from a closed set: this continuity 
is equivalent to \PCA. We also show that if instead of a closed 
set we consider a separably closed set the continuity of the 
function is equivalent to \ACA\ (and hence to the existence of the 
$\inf$ needed for its definition).

The last equivalent condition of the next theorem asserts that the 
sets needed in the straightforward definition of the set of 
$\varepsilon$-witnesses as a $G_\delta$ are indeed open.

\begin{theorem}[\RCA]\label{cont-d}
The following are equivalent:\begin{enumerate}
\item \PCA.
\item For every complete separable metric space $\A$ and every 
closed set $C$ in $\A$ there exists a continuous function $f_C: \A 
\to \R$ such that for every $x \in \A$ we have $f_C(x) = \inf 
\set{d(x,y)}{y \in C}$.
\item For every complete separable metric space $\A$ and every 
open set $U$ in $\A$ the set $\set{(x,\varepsilon) \in \A \times 
\R}{B(x,\varepsilon) \nsubseteq U}$ is open.
\item For every complete separable metric space $\A$, every open 
set $U$ in $\A$ and every $\varepsilon \in \R$ the set $\set{x \in 
\A}{B(x,\varepsilon) \nsubseteq U}$ is open.
\end{enumerate}\end{theorem}
\begin{pf}
$(1) \implies (2)$. Let $\A$ and $C$ be given. Let $\Phi$ be a set 
which enumerates all quadruples $(a,r,c,s) \in A \times \Q^+ 
\times \Q \times \Q^+$ such that
$$\forall b \in A \(d(a,b) < r \imply \exists x \in C\; \(d(b,x) < 
c+s\) \land \forall y \in C \(d(b,y) > c - s\)\).$$
The preceding formula is equivalent to a Boolean combination of 
\P11 formulas and hence $\Phi$ exists within \PCA. Moreover \PCA\ 
shows that $\Phi$ is a code for the function $f_C$.

$(2) \implies (3)$. If $C$ is the complement of $U$ we have that 
$(x,\varepsilon)$ is such that $B(x,\varepsilon) \nsubseteq U$ if 
and only if $f_C(x) < \varepsilon$. If $\Phi$ is a code for the 
continuous function $f_C$ then 
$\set{(x,\varepsilon)}{B(x,\varepsilon) \nsubseteq U} = 
\set{(x,\varepsilon)}{\exists (n,a,r,c,s) \in \Phi\; (d(a,x) < r 
\land c+s < \varepsilon)}$ is open by lemma \ref{s01open}.

$(3) \implies (4)$ is trivial.

$(4) \implies (1)$. We reason within \RCA\ and begin by showing 
that (4) implies \ACA. To this end let $f: \N \to \N$ be a 
one-to-one function: we need to show that the range of $f$ exists. 
$\N$ with the usual metric can be viewed as a complete separable 
metric space $\A$ and we can consider the open set $U = \set{n \in 
\A}{\exists k\; n = f(k)}$: by lemma \ref{s01open} $U$ can be coded 
within \RCA. By (4) let $U'$ be the open set $\set{n \in \A}{B(n,1) 
\nsubseteq U}$. It is immediate to check that for every $n$ we have
$$\exists k\; f(k) = n \text{\quad iff \quad} n \notin U'.$$
The right-hand side of the above equivalence gives a \P01 
definition of the range of $f$, that therefore exists within \RCA\ 
by \D01-comprehension.

Now we can prove within \ACA\ that (4) implies \PCA. It is 
well-known that \PCA\ is equivalent over \RCA\ (and, a fortiori, 
over \ACA) to the statement that if $\<T_n: n \in \N\>$ is an 
infinite sequence of trees of finite sequences of natural numbers 
then the set $X = \set{n}{T_n \text{ is well-founded}}$ exists. Fix 
a sequence  of trees $\<T_n: n \in \N\>$ and let $T = \{\<\>\} \cup 
\set{s \in \seq}{\<s(1), \dots, s(\lh(s)-1)\> \in T_{s(0)}}$. We 
work in the Baire space $\Bai$ (with the metric described in 
section \ref{basic}) and consider the open set $U = \set{x \in 
\Bai}{\exists n\; x[n] \notin T}$, whose elements are all infinite 
sequences which are not paths through $T$. Once more lemma 
\ref{s01open} insures that we can find a code for $U$. By (4) let 
$U'$ be the open set $\set{x \in \Bai}{B(x,1) \nsubseteq U}$. Let 
$x_n \in \Bai$ be the infinite sequence consisting of $n$ followed 
by infinitely many $0$'s. It is immediate to check that for every 
$n$ we have
$$T_n \text{ is well-founded} \text{\quad iff \quad} x_n \notin 
U'.$$
The right-hand side of the above equivalence gives an arithmetical 
definition of $X$, that therefore exists within \ACA.
\end{pf}

\begin{remark}
The open set $U$ used in the second part of the preceding proof is 
the same used by Brown in his proof of theorem \ref{clo->sep}. The 
same open set yields an immediate proof of $(2) \implies (1)$ 
without the need of obtaining \ACA\ first.
\end{remark}

\begin{theorem}[\RCA]\label{cont-d2}
The following are equivalent:\begin{enumerate}
\item \ACA.
\item For every complete separable metric space $\A$ and every 
separably closed set $\cl C$ in $\A$ there exists a continuous 
function $f_{\cl C}: \A \to \R$ such that for every $x \in \A$ we 
have $f_{\cl C}(x) = \inf \set{d(x,y)}{y \in \cl C}$.
\end{enumerate}\end{theorem}
\begin{pf}
$(1) \implies (2)$. Let $\A$ and $C = \<x_n : n \in \N\>$ (the 
code for $\cl C$) be given. A code for the function $f_{\cl C}$ can 
be obtained by setting $(a,r) \Phi (c,s)$ if and only if
$$\exists q \in \Q^+ \forall b \in A \(d(a,b) < r \imply c-s+q < 
\inf \set{d(b,x_n)}{n \in \N} < c+s\).$$
Since \ACA\ proves that every sequence of reals which has a lower 
bound (in this case $0$) has a greatest lower bound and the 
formula is arithmetical we have that the code exists within \ACA.

$(2) \implies (1)$. We use theorem \ref{sep->clo} and prove that 
(2) implies that every separably closed set is closed. This is 
immediate because the complement of $\cl C$ is the preimage by 
$f_{\cl C}$ of the open interval $(0,+\infty)$ and hence it is 
open.
\end{pf}

\begin{remark}
Another proof of $(2) \implies (1)$ of the last theorem consists 
of deducing from (2) the existence of the least upper bound for 
any bounded sequence of real numbers.
\end{remark}


\begin{thebibliography}{99}
\bibitem{ats} {\sc M. Atsuji}, Uniform continuity of continuous 
functions of metric spaces, Pacific J. Math.\ {\bf 8} (1958) 
11--16.

\bibitem{bee1} {\sc G. Beer}, Metric spaces on which continuous 
functions are uniformly continuous and Hausdorff distance, Proc.\ 
Amer.\ Math.\ Soc.\ {\bf 95} (1985) 653--658.

\bibitem{bee2} {\sc G. Beer}, More about metric spaces on which 
continuous functions are uniformly continuous, Bull.\ Austral.\ 
Math.\ Soc.\ {\bf 33} (1986) 397--406.

\bibitem{br-thesis} {\sc D. K. Brown}, {\em Functional Analysis in 
Weak Subsystems of Second Order Arithmetic}, Ph.d.\ thesis, The 
Pennsylvania State University, 1987.

\bibitem{br} {\sc D. K. Brown}, Notions of closed subsets of a 
complete separable metric space in weak subsystems of second order 
arithmetic, in {\em Logic and Computation} (W. Sieg, ed.), 
Contemporary Mathematics 106, American Mathematical Society, 1990, 
pp.\ 39--50.

\bibitem{bs1} {\sc D. K. Brown and S. G. Simpson}, Which set 
existence axioms are needed to prove the separable Hahn-Banach 
theorem?, Ann.\ Pure Appl.\ Logic {\bf 31} (1986) 123--144.

\bibitem{bs2} {\sc D. K. Brown and S. G. Simpson}, The Baire 
category theorem in weak subsytems of second order arithmetic, J. 
Symb.\ Logic {\bf 58} (1993) 557--578.

\bibitem{doc} {\sc M. Do Carmo}, {\em Differential Geometry of 
Curves and Surfaces}, Prenctice Hall, 1976.

\bibitem{hy} {\sc J. G. Hocking and G. S. Young}, {\em Topology}, 
Addison Wesley, 1961.

\bibitem{kel} {\sc J. L. Kelley}, {\em General Topology}, Van 
Nostrand, 1955.

\bibitem{kos} {\sc C. Kosniowski}, {\em A First Course in Algebraic 
Topology}, Cambridge University Press, 1980.

\bibitem{mos} {\sc Y.N. Moschovakis}, {\em Descriptive Set 
Theory\/}, North-Holland, 1980.

\bibitem{st} {\sc N. Shioji and K. Tanaka}, Fixed point theory in 
weak second-order arithmetic, Ann.\ Pure Appl.\ Logic {\bf 47} 
(1990) 167--188.

\bibitem{siZ2} {\sc S. G. Simpson}, Subsystems of $Z_2$ and 
reverse mathematics, appendix to G. Takeuti, {\em Proof Theory}, 
2nd edition, North-Holland, 1986.

\bibitem{sosoa} {\sc S. G. Simpson}, {\em Subsystems of Second 
Order Arithmetic\/}, in preparation.
\end{thebibliography}
\end{document}